\documentclass[mathpazo]{cicp}

\usepackage{amsmath,amsthm,amssymb}
\usepackage{color}
\usepackage{graphicx}
\usepackage{hyperref}
\usepackage{geometry}
\newtheorem{algorithm}{Algorithm}[section]
\numberwithin{equation}{section}
\allowdisplaybreaks 
\newcommand{\comm}[1]{{\color{black}#1}} 
\newcommand{\revise}[1]{{\color{black}#1}} 
\begin{document}
\title{An Algebraic Multigrid Method for Eigenvalue Problems in Some Different Cases}


\author[Zhang N et.~al.]{Ning Zhang\affil{1}, Xiaole Han\affil{2},
       Yunhui He\affil{3},  Hehu Xie\affil{1} and  Chun'guang You\affil{4}\comma\corrauth}
       \address{
       \affilnum{1}\ LSEC, ICMSEC, Academy of Mathematics and Systems Science,
Chinese Academy of Sciences, Beijing 100190, P.R. China,
and School of Mathematical Sciences, University of Chinese Academy of Sciences,
Beijing, 100049, P.R. China \\
\affilnum{2}\ Institute of Applied Physics and Computational Mathematics,
Beijing, 100094, P.R. China \\
           \affilnum{3}\ Department of Mathematics and Statistics,
Memorial University of Newfoundland, St. John's, NL A1C 5S7, Canada\\
\affilnum{4}\ CAEP Software Center for High Performance Numerical Simulation,
Beijing, 100088, P.R. China}
 \emails{{\tt hanxiaole@lsec.cc.ac.cn} (X.~Han), {\tt yunhui.he@mun.ca} (Y.~He),
          {\tt hhxie@lsec.cc.ac.cn} (H.~Xie), {\tt youchg@lsec.cc.ac.cn} (C.~You), {\tt zhangning114@lsec.cc.ac.cn} (N.~Zhang)}


\begin{abstract}
The aim of this paper is to develop an algebraic multigrid method to solve eigenvalue problems based on
the combination of the multilevel correction scheme and the algebraic multigrid method for linear equations.
Our approach uses the algebraic multigrid method setup procedure to construct the hierarchy and the intergrid transfer operators.
In this algebraic multigrid scheme, a large scale eigenvalue problem is solved by some algebraic multigrid smoothing steps in the hierarchy
and very small-dimensional eigenvalue problems. To emphasize the efficiency and flexibility of the proposed method,
here we consider a set of test eigenvalue problems, discretized on unstructured meshes, with different shape of domain, singularity,
and discontinuous parameters. Moreover, global convergence independent of the number of desired eigenvalues is obtained.
\end{abstract}


\maketitle


\section{Introduction}\label{sec;introduction}
Algebraic multigrid (AMG) method was first introduced in \cite{IntroAMG1982}, where the main idea is to design
a similar multigrid method for matrices. Since, however, there is no geometric background, the convergence
has been proved only for some special matrices, such as symmetric positive definite
$M$-matrices with weak diagonal dominance \cite{Ruge1987} and without the assumption of
$M$-matrices in \cite{Huang91ConvergenceWD,Mac2014, Xu2016}. The essential difficulties
for AMG method lie in the choice of coarse grids and intergrid transfer operators,
which fully depend on the understanding of algebraic smooth error under certain smoothing processes.
The classical coarsening strategy was introduced in \cite{Ruge1987}, and others like aggregation and smooth
aggregation in \cite{Aggregation10,SmoothAggregation_Elliptic96}, compatible relaxation \cite{CR_BJFR10,CR_Livne04},
based on element interpolation \cite{AMGe}, energy-based strategy \cite{Energy_AMG06} and so on.
The paper \cite{Robust_AMG} presents some numerical experiments to study the robustness and
scalability of the AMG method. Parallel and adaptive AMG methods have also been studied in \cite{AdaptiveAMG,Boomer}.
So far, due to its simplicity, AMG method has been applied to many problems, such as
\cite{AMG4DG,AMG4Markov,AMG4Helmholtz} and there have been developed many software.

In this paper, we are interested in computing $q$ eigenpairs (maybe not the smallest magnitude)
of the following generalized eigenvalue problem:
Find $(\lambda^{(j)},u^{(j)})\in \mathbb{R}\times \mathbb{R}^N,\ \ j=1,2,\cdots,q$
such that $(u^{(j)})^TMu^{(k)}=\delta_{jk},\ \ j,k=1,2,\cdots,q$ and
\begin{equation}\label{eigenvalue_problem}
Au^{(j)} = \lambda^{(j)} M u^{(j)},\ \ \ \ \ j=1,2,\cdots,q,
\end{equation}
where $A$ is a real symmetric positive definite $N\times N$ matrix,
and $M$ is a real symmetric semi-positive $N\times N$ matrix.
The concerned generalized eigenvalue problem (\ref{eigenvalue_problem}) usually aries
from the discretization of the elliptic partial differential equations involved in several scientific and
theoretical fields such as material sciences, electromagnetics, quantum chemistry,
acoustic and data science etc.. These important applications usually require
high resolution which means the discretization results in large scale algebraic eigenvalue
problems. Thus, the design of efficient eigensolvers with nearly optimal computation complexity
is urgently needed.

\comm{It is a natural idea to use AMG and MG methods for eigenvalue problems
\cite{BrandtMcCormickRuge,CaiMandelMcCormick,GongHanSunZhang,Hackbusch,Shaidurov,XieXieYinYue,ZhaiXieZhangZhang}.
A very good review of the
application of AMG methods to eigenvalue problems is given in \cite{KnyazevNeymeyr} and references cited therein.
Roughly speaking, in these methods, the AMG strategy is adopted as only the solver for linear equations
in the inner iterations combined with some types of outer iterations,
such as inverse power, shift-and-inverse, Rayleigh-quotient,
locally optimal block preconditioned conjugate gradient and so on,
for eigenvalue problems.
But the application of the AMG method does not lead to new eigensolver (outer iteration).}
Recently, a type of multilevel correction method
\cite{KuzelVanek,LinXie_2011,LinXie_2012,LinXie_MultiLevel,Pultarova,VanekPultarov,Xie_JCP,Xie_IMA} is proposed to solve eigenvalue problems.
This method is based on a new understanding of Aubin-Nitsche technique in the finite element method \cite{LinXie_2011}.
In contrast to the methods reviewed in \cite{KnyazevNeymeyr} for which AMG is only used as a preconditioner
of the stiffness matrix, the coarse space from the multigrid method plays the role to enriching
the working subspace for the eigenvalue problem solving \cite{HeHongXieYouYue}.
Then in this multilevel correction scheme, the solution of eigenvalue problem on the finest level mesh
can be reduced to a sequence of solutions of standard boundary value problems on the multilevel meshes
and some solutions of the eigenvalue problems on a very low dimensional space.
Therefore, the computational work and required memory can arrive at the optimality.
\comm{The above description shows that the application of the multigrid method
in the multilevel correction scheme leads to a new eigensolver.}


Motivated by the AMG method for boundary value problems and the multilevel correction method,
we design a type of AMG method for eigenvalue problems. The proposed AMG method in this paper
can compute several eigenpairs, \comm{which may be not the smallest magnitude},
and permits the free choice of the eigensolvers for the included low dimensional eigenvalue problems.
With simple Gauss-Seidel relaxation as a smoothing process, the AMG still achieves robust convergence.
Our AMG method also allows us to add simple strategies so that problems descritized on completely unstructured grids
can also be solved efficiently. Without the geometric background, this AMG method has a wide range of applicability.
Another aim of this paper is to investigate the efficiency of the AMG method for eigenvalue problems.
We test several different types of eigenvalue problems. Numerical results show that the time consumption and
iteration numbers are almost optimal, and the convergence is independent of the number of the computed eigenvalues.

The rest of this paper is organized as follows. In the next section, we introduce
the classical AMG method, mainly the constructing of ``coarse-grid''.
An AMG algorithm for solving the eigenvalue problem is proposed
and analyzed in Section \ref{AMG-eig-alg}. In Section \ref{sec-numeri}, some numerical tests are presented to
validate the efficiency of the proposed algorithm.
Some concluding remarks are given in the last section.

\section{Algebraic multigrid hierarchy}\label{sec;preliminary}
This section is devoted to introducing the classical AMG method which
aims at solving the ill-conditioned linear system $Au=f$ similarly to the
geometric multigrid (GMG) method.

\subsection{Standard coarsening}
Since there is no real geometric background, the main content is to
determine the ``coarse-grid" directly from the matrix $A = (a_{ij})$.
By analogy, we define grid points, $\Omega$, as the indices $\{1,2,\cdots,N\}$ of
$u=(u_1,u_2,\cdots,u_N)^T$, and choose a subset of $\Omega$ as the coarse
grid points according to the undirected adjacency graph of the matrix $A$.

In order to derive a coarse level system, $\Omega$ is needed to split into
two disjoint subsets $\Omega = C \cup F$ with $C$ representing the coarse grid points
and $F$ being the complementary.
Following \cite{Ruge1987,Stuben2001Review,Stuben2001Introduction}, we define the strong dependent set
\begin{equation}
S_i:=\Big\{j\ \ \big|\ -a_{ij}\ge\theta\max\limits_{\substack{a_{i\ell}< 0\\ \ell\neq i}}|a_{i\ell}|\Big\},
\end{equation}
and the strong influence set
$S_i^T:=\{j\big|i\in S_j\}$ with fixed $0<\theta<1$ (usually $0.25$).
The measure of how valuable a variable $i$ is as a coarse grid point
can be represented by $|S_i^T|$, i.e., the the number of elements that  $S_i^T$ contains.
\comm{About the detailed coarsening process, please refer to \cite{Ruge1987,Stuben2001Review,Stuben2001Introduction}.}
%

\subsection{Interpolation}
After the coarse variables set $C$ and its complement $F$ have been constructed,
we can define interpolation from coarse level (with meshsize $H$)
to fine level (with meshsize $h$). For any vector in the coarse grid, $e^H$,
the interpolation (prolongation) operator to fine grid can be defined as follows:
\begin{equation}\label{Interpolation_Def}
\left(I_H^h e^H\right)_i=\left\{
\begin{array}{rc}
e^H_i,                           & \text{if}\ \  i\in C,\vspace{0.5em}\\
\sum\limits_{j\in P_i}\omega_{ij}e^H_j, & \text{if}\ \  i\in F,
\end{array}
\right.
\end{equation}
where $P_i$ is some small sets of interpolation points $C\cap N_i$, where
$N_i:=\{j\in\Omega: j\neq i, a_{ij}\neq 0\}$ denotes the neighborhood of point $i$, $\forall \  i\in \Omega$ .

If choosing the simplest case $P_i = S_i\cap C$, direct interpolation is applied immediately,
as described in \cite{Stuben2001Introduction}. More precisely, the interpolation weights are
\begin{equation}\label{Interpolation_Dir}
\omega_{ij}=\left\{
\begin{array}{rc}
-\alpha_i \displaystyle{\frac{a_{ij}}{a_{ii}}}, & \text{if}\ \  a_{ij}<0, \vspace{0.5em}\\
-\beta_i  \displaystyle{\frac{a_{ij}}{a_{ii}}}, & \text{if}\ \  a_{ij}>0,
\end{array}
\right.
\end{equation}
with
\begin{equation}
\alpha_i=\displaystyle{\frac{\sum\limits_{j\in N_i, a_{ij}<0}a_{ij}}{\sum\limits_{\ell\in P_i, a_{i\ell}<0}a_{i\ell}}}
\qquad\text{and}\qquad
\beta_i=\displaystyle{\frac{\sum\limits_{j\in N_i, a_{ij}>0}a_{ij}}{\sum\limits_{\ell\in P_i, a_{i\ell}>0}a_{i\ell}}}.
\end{equation}
The leading coefficients $\alpha_i$ and $\beta_i$ are chosen to ensure that the approximation
interpolates constants exactly for the case of zero row sum matrices.
For more details and other interpolation types such as standard interpolation,
please see references \cite{Ruge1987,Stuben2001Review,Stuben2001Introduction}.

\begin{remark}
Of course, different type of coarsening and interpolation strategies can also be adopted for
different types of matrices.
\end{remark}

\subsection{Coarse problem}
The AMG setup procedure finally constructs a hierarchy of vector spaces which are
indexed by $k=1,2,\cdots,n$, where $k=1$ denotes the finest level and $k=n$ denotes the coarsest level.
Denote the finest grid $\Omega_1=\Omega$ and corresponding matrices $A_1=A$, $M_1=M$. Based on $A_1$,
the AMG scheme builds up the prolongation and restriction operators $I_{k+1}^k$ and
$I_k^{k+1}:=(I_{k+1}^k)^T$, respectively for $k=1,2,\cdots,n-1$. The coarse matrices are defined with the
Galerkin projection as follows:
\begin{equation}\label{Coarse_Matrices}
A_{k+1}=I_{k}^{k+1}A_kI_{k+1}^k,\ \ {\rm and}\ \
M_{k+1}=I_{k}^{k+1}M_kI_{k+1}^k,\ \  {\rm for}\ \ k=1,2,\cdots,n-1.
\end{equation}
In the following, $d_1,d_2,\cdots,d_n$ denote the dimensions of problems defined on grids $\Omega_1,\Omega_2,\cdots,\Omega_n$, respectively.

\section{AMG algorithm for eigenvalue problems}\label{AMG-eig-alg}
In this section, we propose an AMG method for solving eigenvalue problems. Similar to the geometric
case in \cite{Xie_JCP}, assume we have obtained eigenpair approximations
$\{\lambda_k^{(j,\ell)},u_k^{(j,\ell)}\}_{j=1}^{q}$ to our desired eigenpairs.
Now we introduce an AMG correction step to improve their accuracy.
\begin{algorithm}\label{AMG_One_Correction_Step_Alg}
AMG Correction Step
\begin{enumerate}
\item For $j=1,\cdots,q$ Do\\
\ \ \ Solve the following linear equation by AMG iterations
\begin{eqnarray}\label{Linear_Equation}
\ \ \ A_k\widehat{u}_{k}^{(j,\ell+1)}=\lambda_k^{(j,\ell)}M_ku_k^{(j,\ell)}.
\end{eqnarray}
Perform $m$ AMG iteration steps with the initial value $u_k^{(j,\ell)}$ to obtain
a new eigenfunction approximation $\widetilde{u}_k^{(j,\ell+1)}$ which is denoted by
\begin{eqnarray}\label{inexact_inverse_iter}
\widetilde{u}_k^{(j,\ell+1)}={\rm AMG}(k,\lambda_k^{(j,\ell)}u_k^{(j,\ell)},u_k^{(j,\ell)},m),
\end{eqnarray}
where $k$ denotes the working level $\Omega_k$ for the AMG iteration,
$\lambda_k^{(j,\ell)}u_k^{(j,\ell)}$ leads to the right hand side term of the linear equation,
$u_k^{(j,\ell)}$ denotes the initial guess and $m$ is the number of AMG cycles.

\item Set $V_{k,\ell+1}=[\widetilde{u}_k^{(1,\ell+1)},\cdots,\widetilde{u}_k^{(q,\ell+1)}]$ and
construct two matrices $A_{n,k}^{(\ell+1)}$ and $M_{n,k}^{(\ell+1)}$ as follows
\begin{equation}
A_{n,k}^{(\ell+1)} = \left( \begin{array}{cc}
A_n & I_k^nA_kV_{k,\ell+1} \\
V_{k,\ell+1}^TA_kI_n^k & V_{k,\ell+1}^TA_kV_{k,\ell+1}
\end{array}   \right),
\end{equation}
and
\begin{equation}
M_{n,k}^{(\ell+1)} = \left( \begin{array}{cc}
M_n & I_k^nM_kV_{k,\ell+1} \\
V_{k,\ell+1}^TM_kI_n^k & V_{k,\ell+1}^TM_kV_{k,\ell+1}
\end{array}
\right).
\end{equation}
Solve the following eigenvalue problems: Find $\{\lambda_k^{(j,\ell+1)},x_k^{(j,\ell+1)}\}_{j=1}^q$
such that\\ $(x_k^{(j,\ell+1)})^TM_{n,k}^{(\ell+1)}x_k^{(j,\ell+1)}=1$ and
\begin{eqnarray}\label{Eigenvalue_Problem_Coarse}
A_{n,k}^{(\ell+1)}x_k^{(j,\ell+1)} = \lambda_k^{(j,\ell+1)}M_{n,k}^{(\ell+1)}x_k^{(j,\ell+1)},\ \ \ j=1,\cdots,q.
\end{eqnarray}
Select the desired eigenpairs $\{\lambda_k^{(j,\ell+1)},x_k^{(j,\ell+1)}\}_{j=1}^q$
and do the following computation:\\
For $j=1,\cdots,q$ Do
\begin{eqnarray*}
u_k^{(j,\ell+1)} = I_n^{k}x_{k}^{(j,\ell+1)}(1:d_n) + V_{k,\ell+1}x_k^{(j,\ell+1)}(d_n+1:d_n+q).
\end{eqnarray*}
\end{enumerate}
Summarize the above two steps by defining
\begin{equation*}
\{\lambda_k^{(j,\ell+1)},u_k^{(j,\ell+1)}\}_{j=1}^q =
{\rm AMGCorrection}\big(n,k,\{\lambda_k^{(j,\ell)},u_k^{(j,\ell)}\}_{j=1}^q\big).
\end{equation*}
\end{algorithm}
\comm{Due to the construction of (\ref{Linear_Equation}), the result
$\widetilde{u}_k^{(j,\ell+1)}$ can be viewed as performing a few inexact inverse power
iterations on the given approximation $u_k^{(j,\ell)}$.
Denote $V_n$ the coarsest space of AMG.
The matrices $A_{n,k}^{(\ell+1)}$ and $M_{n,k}^{(\ell+1)}$ are exactly
the Galerkin projections on the augmented subspace
$V_n + \text{span}\{V_{k,\ell+1}\}$.
Thus, the results of Algorithm \ref{AMG_One_Correction_Step_Alg} will never be less accurate
than the (block) inverse power method, where the projection subspace
is $\text{span}\{V_{k,\ell+1}\}$.
Moreover, the approximate low frequency information contained in $V_n$ brings
an improvement of the convergence rate of Algorithm \ref{AMG_One_Correction_Step_Alg} rather than
the usual (block) inverse power method (cf. \cite{LinXie_2011,LinXie_2012,LinXie_MultiLevel,Xie_JCP,Xie_IMA}).
The efficiency and convergence rate
of Algorithm \ref{AMG_One_Correction_Step_Alg} will be shown in section \ref{sec-numeri}.}


\begin{remark}
In solving of (\ref{Linear_Equation}), we may combine other efficient iterate strategies,
such as shifting and polynomial filtering, to obatin better results.
More general correction subspaces can also be considered, for example,
  $V_n + \text{KrylovSpaces}$.
\end{remark}

Based on the above algorithm, we can construct an AMG method for eigenvalue problems
which is a combination of the nested technique and the AMG correction step defined
by Algorithm  \ref{AMG_One_Correction_Step_Alg}.
\begin{algorithm}\label{AMG_Eigenvalue_Multi_Alg}
AMG Eigenvalue Solver
\begin{enumerate}
\item Define the following small dimensional eigenvalue problem
in the $n_1$-th grid $\Omega_{n_1}$ ($n_1\leq n$):
Find $\{\lambda_{n_1}^{(j)},u_{n_1}^{(j)}\}_{j=1}^{q}$ such that
$(u_{n_1}^{(j)})^TM_{n_1}u_{n_1}^{(j)}=1$ and
\begin{eqnarray}
A_{n_1} u_{n_1}^{(j)} = \lambda_{n_1}^{(j)} M_{n_1} u_{n_1}^{(j)}.
\end{eqnarray}
Solve these eigenvalue problems to get eigenpair approximations
$\{\lambda_{n_1}^{(j)},u_{n_1}^{(j)}\}_{j=1}^q$ which are approximations to our desired eigenpairs.
\item For $k=n_1-1,\cdots,1$, perform the following correction steps
\begin{itemize}
\item [(a)] Set $\lambda_k^{(j,0)}=\lambda_{k+1}^{(j)}$ and $u_k^{(j,0)}=I_{k+1}^ku_{k+1}^{(j)}$
for $j=1,\cdots,q$.
\item [(b)] Do the following correction iteration for $\ell=0,\cdots,p_k-1$
\begin{equation*}
\{\lambda_k^{(j,\ell+1)},u_k^{(j,\ell+1)}\}_{j=1}^q =
{\rm AMGCorrection}\big(n,k,\{\lambda_k^{(j,\ell)},u_k^{(j,\ell)}\}_{j=1}^q\big).
\end{equation*}
\item [(c)] Set $\lambda_k^{(j)}=\lambda_k^{(j,p_k)}$ and $u_k^{(j)}=u_k^{(j,p_k)}$ for
$j=1,\cdots,q$.
\end{itemize}
\end{enumerate}
Finally, we obtain eigenpair approximations  $\{\lambda_1^{(j,\ell+1)},u_1^{(j,\ell+1)}\}_{j=1}^q$
in the finest level grid $\Omega_1$.
\end{algorithm}
Different from the GMG method \cite{Xie_JCP}, we do not have the exact
prolongation and restriction operators. In the practical computation,
we can choose suitable iterations $p_k$ to meet the accuracy requirement.
Compared with other AMG methods, the proposed method here only needs to do smoothing iterations for
the standard elliptic type of linear equations and the AMG method can act as a black-box.
The involved small scale eigenvalue problems can be solved by any eigensolver which can also
act as a black-box in our method. Furthermore, different form the classical eigensolvers
such as Lanczos and Arnoldi, the required memory for the eigenpair solving is only about $qN$.
Inspired by the analysis for GMG method \cite{Xie_JCP,Xie_IMA}, it is well known that
 AMG method can have a very good convergence rate if the coarse grids capture
the low frequency information of the finest grid well.

\section{Numerical results}\label{sec-numeri}
In this section, we illustrate our method by six different examples.
In this section, we present six numerical examples to investigate the efficiency of the proposed algorithm.
These examples include the algebraic eigenvalue problems from the discretization of the Poisson eigenvalue problem defined on
different domains by the linear finite element method. The AMG method will be tested on the structure and unstructure meshes.
Here, we use {\it globally uniform convergence} to denote that the numerical method has the same convergence rate for different
number of computed eigenpairs.


\subsection{Default implementation settings}
Some parameters have to be set to finish the AMG setup phase and the eigenvalue algorithm.
Recall that we want to compute the first $q$ eigenpairs of  generalized eigenvalue problems.
The main default settings and procedures are listed below.

\begin{itemize}
\item Strong dependent/influence set threshold: $\theta = 0.25$.
\item Interpolation type: direct interpolation.
\item Linear solver: 1 AMG V-cycle:
\begin{itemize}
\item Pre-smoothing: 1 Gauss-Seidel iteration.
\item Post-smoothing: 1 Gauss-Seidel iteration.
\end{itemize}
\item The initial eigen solver level: $n_1 = n$.
\item Correction number on each level: $p_{n_1-1} = p_{n_1-2} = \cdots = p_2 = 1$, $p_1 \leq 20$.
\item Direct eigen solver: Arnoldi process by the ARPACK library. The $j$-th eigenvalue on the finest
level is denoted by $\lambda_{j}^{\text{dir}}$.
\item $q$ is the number of desired eigenvalues, ranging from $1$ to $30$.
\item Total errors on the finest level between the eigenvalue approximations by the AMG method and direct solver:
$e_{\ell} = \sum\limits_{j=1}^{q}|\lambda_j^{\ell}-\lambda_j^{\text{dir}}|$
where $\lambda_j^{\ell} = \lambda_1^{(j,\ell)}$, $\ell = 1,2,\cdots,p_1$.
\item Total error tolerance: $\tau \leq 10^{-9}$.
\item Average convergence ratio: ${\rm ratio}=\left(\frac{\displaystyle{e_{p_1}}}{\displaystyle{e_1}}\right)^{\frac{1}{p_1-1}}$.
\item The stiffness matrix $A$ and mass matrices $M$ are obtained
by discretizing the following problems respectively with the linear finite element method \cite{Ciarlet1978}.
\item Test machine: Intel Xeon E5-2620 2.00GHz,
two 6-core dual thread CPUs, 72G memory.
\end{itemize}

\subsection{Example 1: model eigenvalue problem}
The first mode problem is the most elementary eigenvalue problem:
Find $(\lambda,u)$ such that
\begin{equation}\label{model problem}
\left\{
\begin{array}{rrrlr}
-\Delta u                  &=& \lambda u,\ &{\rm in} &\Omega,\\
u                  &=&         0,\ &{\rm on} &\partial\Omega,\\
\int_{\Omega} u^2 d\Omega &=& 1,\ &         &
\end{array}
\right.
\end{equation}
where $\Omega=(0,1)\times(0,1)$. The mesh is generated by uniform refinement beginning with mesh size $h=1$.
The dimensions on each level of grids
are [4198401, 2095105, 525309, 131581, 33021, 8569, 2109, 542] and the total level $n=8$.

\begin{figure}[http]
\begin{center}
\includegraphics[width=7cm,height=6cm]{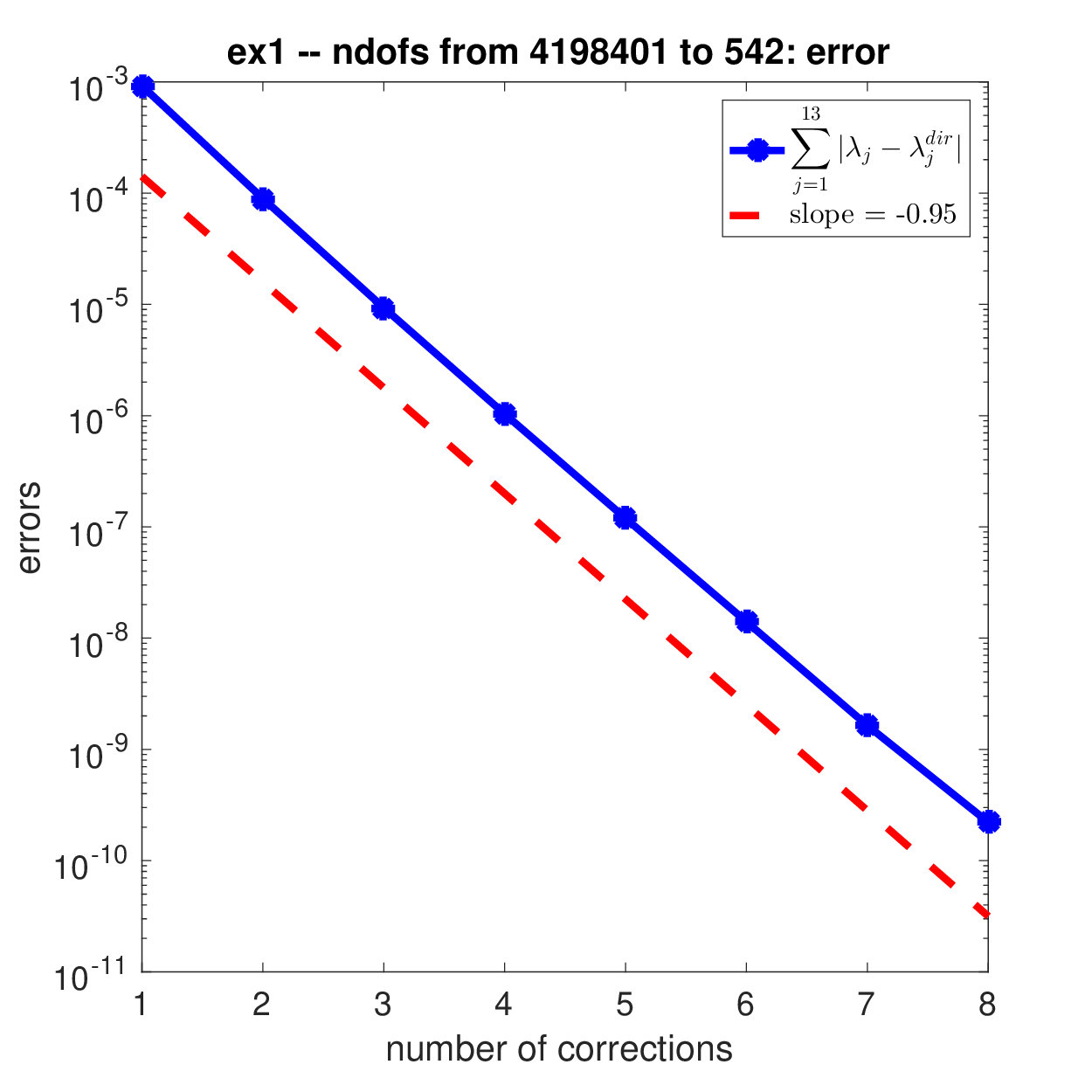}
\includegraphics[width=7cm,height=6cm]{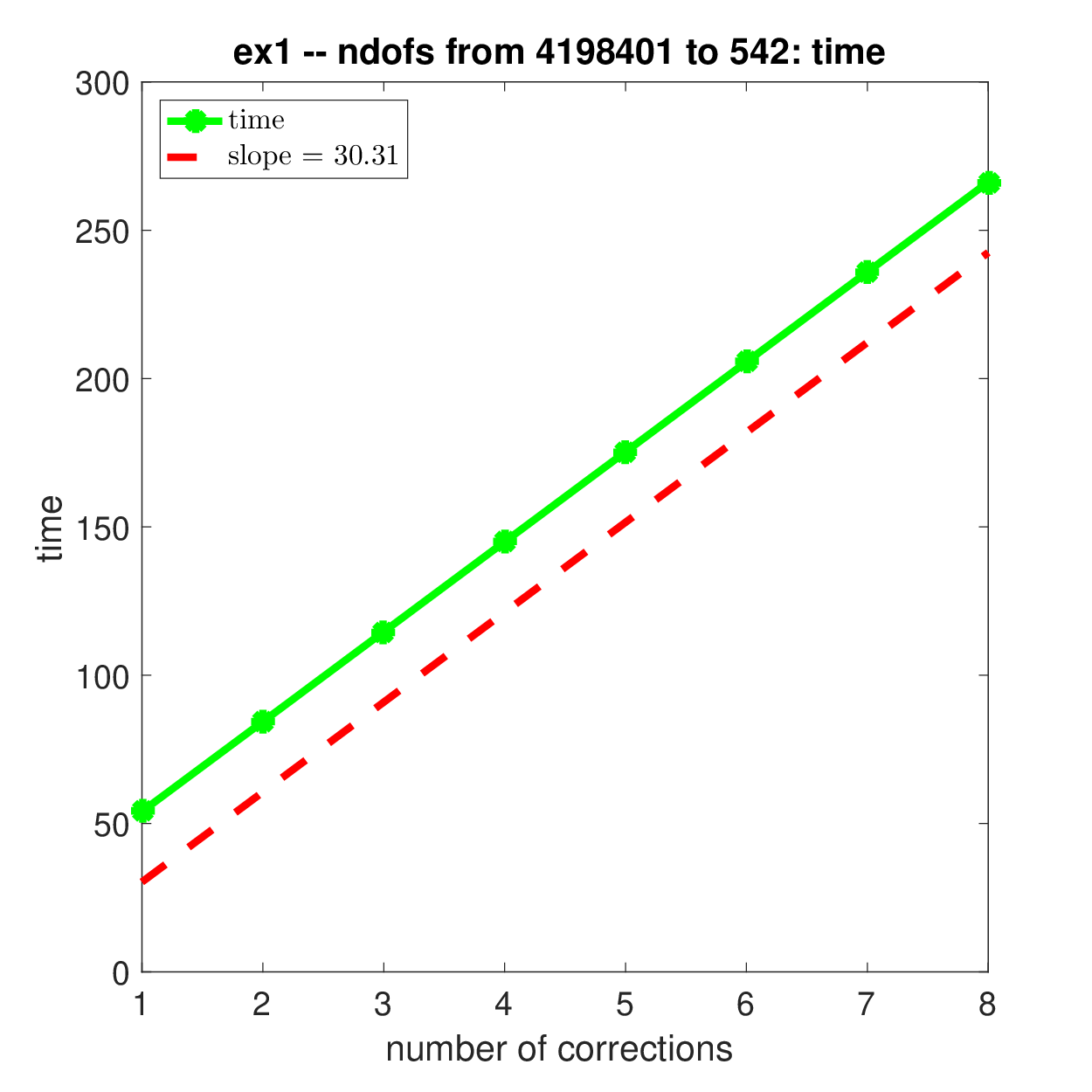}
\caption{Problem (\ref{model problem}), example 1 -- The algebraic errors and CPU time (in second)
of the AMG method for the first $13$ eigenvalues
on the uniform refinement mesh.}\label{example 1 figure}
\end{center}
\end{figure}

Table \ref{example 1 table} shows that the average convergence ratio, total errors, and iteration numbers for different desired amount
of eigenvalues. Different with Krylov methods where more than $q$ eigenvalues are actually computed,  Algorithm \ref{AMG_Eigenvalue_Multi_Alg} only needs to compute the first $q$ eigenvalues. Generally, the more closer $\lambda_q$ and $\lambda_{q+1}$ is, the more difficultly it performs when solving the first $q$ eigenvalues.
However, from Table \ref{example 1 table}, the overall behavior is similar--we get the {\it globally uniform convergence} ratio,
i.e., about $0.11$ for different number of computed eigenpairs. 
To see the optimal \revise{computational} complexity, Figure \ref{example 1 figure} gives the algebraic errors and CPU time (in seconds)
for the first $13$ eigenvalues. From the left of Figure \ref{example 1 figure},  we see the error of the eigenvalues is reduced linearly with the iteration numbers. The right of Figure \ref{example 1 figure} shows that the CPU time grows perfectly linearly with the number of iteration
which is faster than the direct method which needs $1737.56$ seconds under the same error tolerance.

\begin{table}[http]
\begin{center}
\begin{tabular}{c|c|c|c|c|c}
\hline
$q$&       $\lambda_q$     & $\lambda_q / \lambda_{q+1}$ &    ratio   & iter & total error   \\ \hline
1  &   19.739351152446577  &     0.399999830658111       &  0.110359  &   6  &  0.19e-09     \\
2  &   49.348398772994052  &                             &            &      &               \\
3  &   49.348426661655587  &     0.624999514802597       &  0.109683  &   7  &  0.12e-09     \\
4  &   78.957543954641395  &     0.800000031375434       &  0.110713  &   7  &  0.21e-09     \\
5  &   98.696926072478092  &                             &            &      &               \\
6  &   98.696926072787591  &     0.769230457062699       &  0.111902  &   7  &  0.45e-09     \\
7  &  128.306055963593963  &                             &            &      &               \\
8  &  128.306290898228610  &     0.764706565466666       &  0.110547  &   7  &  0.69e-09     \\
9  &  167.784999753374791  &                             &            &      &               \\
10  &  167.785014924843836  &     0.944442983816993       &  0.112142  &   8  &  0.12e-09    \\
11  &  177.654996436879685  &     0.900000105867284       &  0.115364  &   8  &  0.18e-09    \\
12  &  197.394417265854258  &                             &            &      &              \\
13  &  197.394417273111628  &     0.799999907230014       &  0.113346  &   8  &  0.22e-09    \\
14  &  246.743050204326266  &                             &            &      &              \\
15  &  246.743972243036808  &     0.961541879013450       &  0.112188  &   8  &  0.27e-09    \\
16  &  256.612819086151660  &                             &            &      &              \\
17  &  256.612819086431614  &     0.896550451610119       &  0.113451  &   8  &  0.37e-09    \\
18  &  286.222396771513957  &                             &            &      &              \\
19  &  286.222479299660563  &     0.906247979981473       &  0.112243  &   8  &  0.42e-09    \\
20  &  315.832405282174534  &     0.941176611569284       &  0.111708  &   8  &  0.45e-09    \\
21  &  335.571880346205262  &                             &            &      &              \\
22  &  335.571880413890710  &     0.918920429907456       &  0.115726  &   8  &  0.69e-09    \\
23  &  365.180563509384399  &                             &            &      &              \\
24  &  365.180577155566198  &     0.924998520135014       &  0.110632  &   8  &  0.59e-09    \\
25  &  394.790444748240247  &                             &            &      &              \\
26  &  394.790444752357189  &     0.975610272573643       &  0.114287  &   8  &  0.87e-09    \\
27  &  404.659991649028711  &                             &            &      &              \\
28  &  404.662532642826989  &     0.911114427716026       &  0.113430  &   8  &  0.95e-09    \\
29  &  444.140187371669356  &                             &            &      &              \\
30  &  444.140438372226185  &     0.900001567378938       &  0.138346  &   9  &  0.60e-09    \\ \hline
\end{tabular}
\caption{Problem (\ref{model problem}), example 1 -- Results about the algebraic errors on the unit
square with uniform refinement mesh.}\label{example 1 table}
\end{center}
\end{table}

\subsection{Example 2: Poisson eigenvalue problem on L-shaped domain}
In this example, we test the performance of Algorithm \ref{AMG_Eigenvalue_Multi_Alg}
for the problem (\ref{model problem}),  but on an L-shaped domain
$\Omega=[-1,1]\times [-1,1]\backslash(0,1)\times(-1,0)$.
The mesh is generated by uniform refinement beginning with mesh size $h=1$.
The dimensions on each level of grids
are [3149825, 1570817, 394236, 98812, 24828, 6518, 1596] and the total level $n=7$.

\begin{figure}[http]
\begin{center}
\includegraphics[width=7cm,height=6cm]{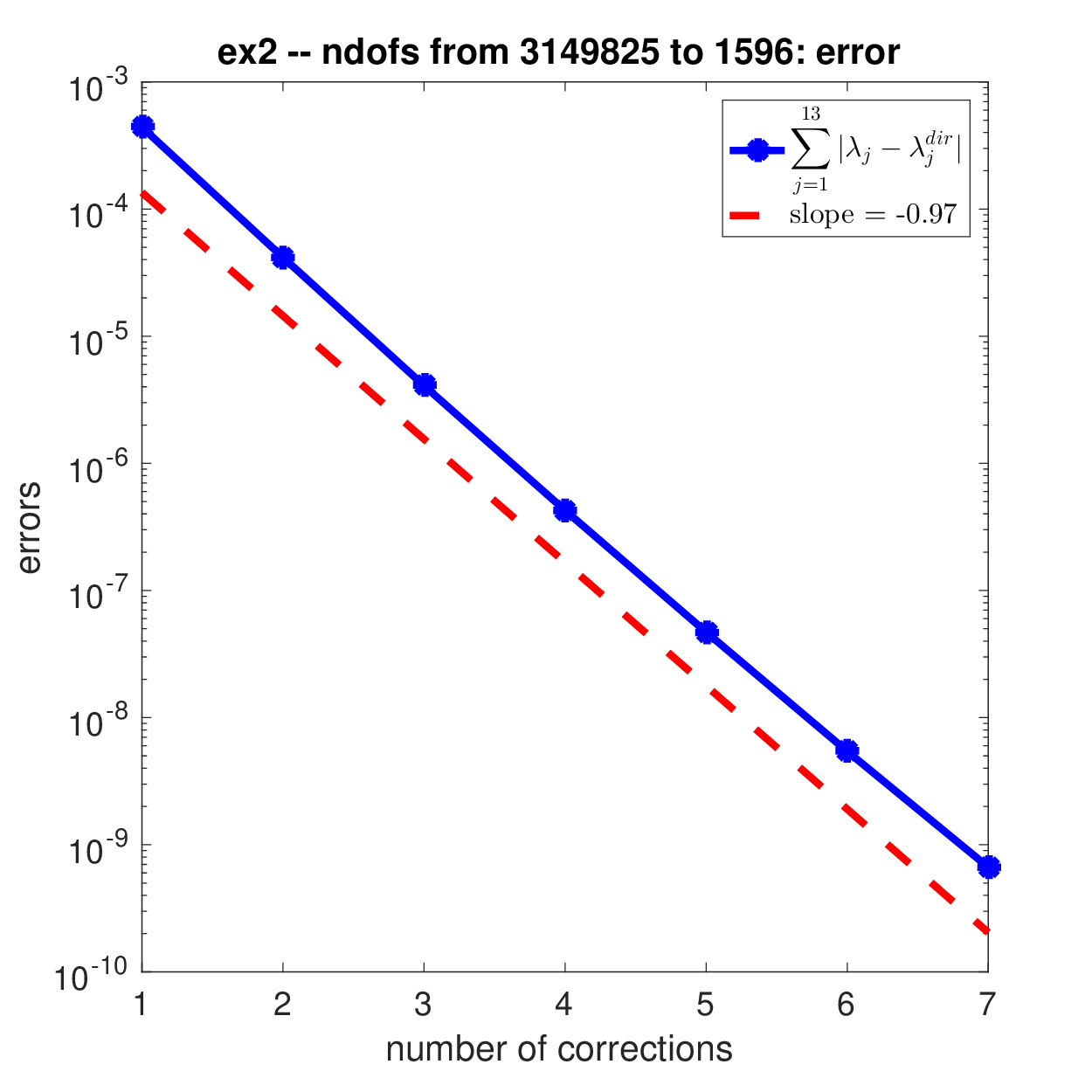}
\includegraphics[width=7cm,height=6cm]{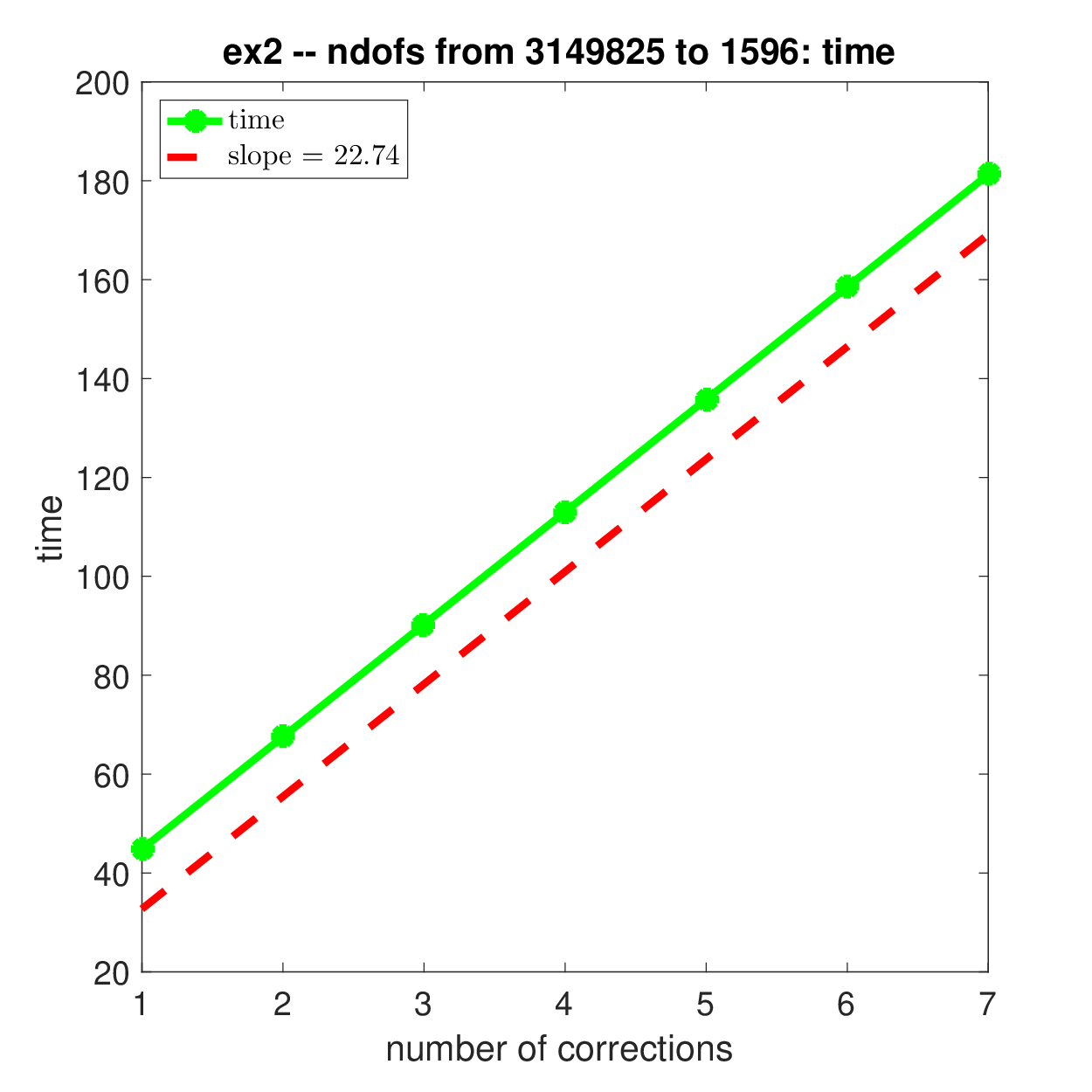}
\caption{Problem (\ref{model problem}), example 2 -- The algebraic errors and CPU time (in second)
of the AMG method for the first $13$ eigenvalues
on the uniform refinement mesh.}\label{example 2 figure}
\end{center}
\end{figure}

In Table \ref{example 2 table}, numerical results are also listed for different desired eigenvalues using
our AMG algorithm. Here, again, we see {\it globally uniform convergence}, but
it is worth noting that the dimension on the coarsest level are increased to $1596$.
It is reasonable due to the L-shaped domain. However, the  behaviors of algebraic errors
and CPU time for the first $13$ eigenvalues are similar to Example 1. Again, the timings are close to
scaling linearly with the iteration numbers, see Figure \ref{example 2 figure}. For comparison, the direct method need
936,94 seconds for the first $13$ eigenpairs under the same error tolerance.

\begin{table}[http]
\begin{center}
\begin{tabular}{c|c|c|c|c|c}
\hline
$q$ &       $\lambda_q$     & $\lambda_q / \lambda_{q+1}$ &    ratio   & iter & total error  \\ \hline
1  &    9.639941155614530	&   0.634320854275592         &  0.105045  &   6  &  0.12e-09    \\
2  &   15.197263483672709	&   0.769901102357829         &  0.104864  &   6  &  0.28e-09    \\
3  &   19.739241101397241	&   0.668638429761318         &  0.105088  &   6  &  0.48e-09    \\
4  &   29.521547405588866	&   0.925055808457296         &  0.104818  &   6  &  0.73e-09    \\
5  &   31.913260946733118	&   0.769457056102703         &  0.107524  &   7  &  0.13e-09    \\
6  &   41.475038397040180	&   0.922721052755031         &  0.107113  &   7  &  0.18e-09    \\
7  &   44.948620466830491	&   0.910845766545385         &  0.106877  &   7  &  0.24e-09    \\
8  &   49.348223505840721	&                             &            &      &              \\
9  &   49.348260690989306	&   0.870179621581926         &  0.106558  &   7  &  0.35e-09    \\
10  &   56.710430199775971	&   0.867439196741824         &  0.106295  &   7  &  0.41e-09    \\
11  &   65.376836108842241	&   0.920039372342096         &  0.105972  &   7  &  0.48e-09    \\
12  &   71.058737347746145	&   0.992811733747242         &  0.105879  &   7  &  0.54e-09    \\
13  &   71.573224743772840	&   0.906477598649387         &  0.107138  &   7  &  0.66e-09    \\
14  &   78.957521785881866	&   0.884144399022638         &  0.106518  &   7  &  0.73e-09    \\
15  &   89.303876010710553	&   0.967460206330204         &  0.106860  &   7  &  0.87e-09    \\
16  &   92.307544461658438	&   0.947893993499531         &  0.106015  &   7  &  0.95e-09    \\
17  &   97.381716831929808	&   0.986674656579670         &  0.189571  &  10  &  0.28e-09    \\
18  &   98.696886742288171	&                             &            &      &              \\
19  &   98.696886759985802	&   0.971356701426527         &  0.108106  &   8  &  0.15e-09    \\
20  &  101.607253663911848	&   0.904219670345751         &  0.107446  &   8  &  0.16e-09    \\
21  &  112.370098767105745	&   0.972722811055895         &  0.110938  &   8  &  0.22e-09    \\
22  &  115.521192152497733	&   0.900353570057789         &  0.112182  &   8  &  0.26e-09    \\
23  &  128.306474249980454	&                             &            &      &              \\
24  &  128.306787514516458	&   0.986060874784479         &  0.111082  &   8  &  0.29e-09    \\
25  &  130.120554212801750	&   0.998942245349293         &  0.108884  &   8  &  0.27e-09    \\
26  &  130.258335573047447	&   0.914382448573862         &  0.109380  &   8  &  0.29e-09    \\
27  &  142.454982350337076	&   0.942648379890289         &  0.109414  &   8  &  0.32e-09    \\
28  &  151.122078379763281	&   0.978332935190765         &  0.109511  &   8  &  0.34e-09    \\
29  &  154.468967509814064	&   0.952336564402184         &  0.108736  &   8  &  0.35e-09    \\
30  &  162.199975600831635	&   0.984711949246173         &  0.112557  &   8  &  0.48e-09    \\ \hline
\end{tabular}
\caption{Problem (\ref{model problem}), example 2 -- Results about the algebraic errors on an
L-shape domain with uniform refinement mesh.}\label{example 2 table}
\end{center}
\end{table}

\subsection{Example 3: Poisson eigenvalue problem with discontinuous parameters I}
Our third example involves discontinuous parameters.  We are concerned with the following Poisson eigenvalue problem on the unit square:
Find $(\lambda, u )$ such that
\begin{equation}\label{poisson problem}
\left\{
\begin{array}{rrrll}
-\nabla\cdot(K\nabla u) &=& \lambda u, \ &{\rm in} &\Omega,\\
             u  &=& 0,         \ &{\rm on} &\partial\Omega,\\
\int_{\Omega} u^2 d\Omega &=& 1,\ &         &
\end{array}
\right.
\end{equation}
where $\Omega=(-1,1)\times (-1,1)$ and the coefficient matrix $K$ has discontinuous elements.
We select two different coefficients $K$ (in this and next subsection) to test the performance of our multigrid algorithm.

Here, we set coefficients $K$ as follows:
\begin{equation*}
\left\{
\begin{array}{ccll}
K &=& \begin{pmatrix}    1  &  0 \\  0 &     1  \end{pmatrix} \ &{\rm in}\ (-1,0) \times ( 0,1) \cup (0,1)\times (-1,0),\\
K &=& \begin{pmatrix} 1000  &  0 \\  0 &  1000  \end{pmatrix} \ &{\rm in}\ [ 0,1] \times [ 0,1],\\
K &=& \begin{pmatrix} 0.001 &  0 \\  0 &  0.001 \end{pmatrix} \ &{\rm in}\ [-1,0] \times [-1,0].
\end{array}
\right.
\end{equation*}
The mesh is generated by uniform refinement beginning with mesh size $h=1$.
The dimensions on each level of grids
are [4198401, 2095105, 525310, 131585, 33024, 8564, 2137] and the total level $n=7$.

\begin{figure}[http]
\begin{center}
\includegraphics[width=7cm,height=6cm]{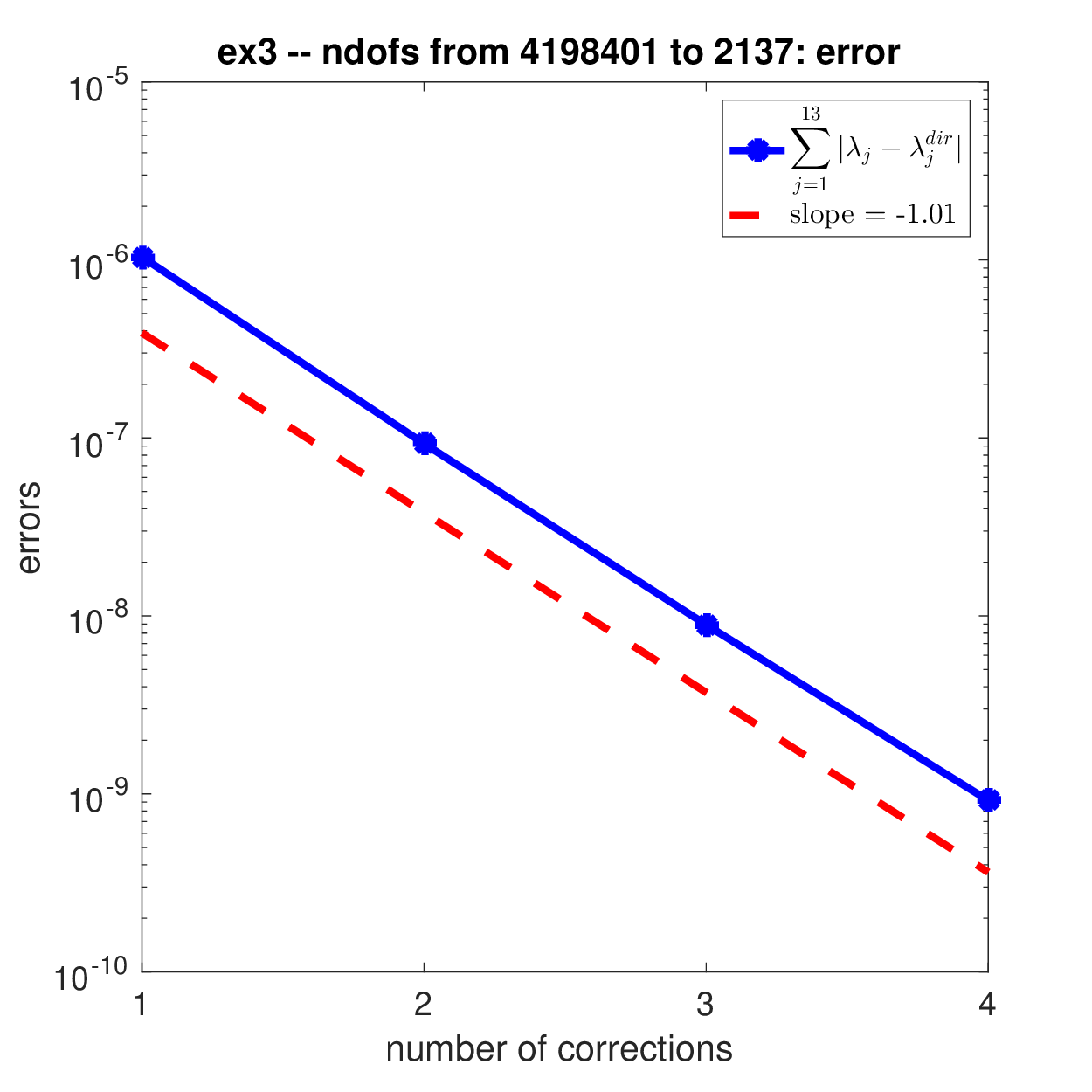}
\includegraphics[width=7cm,height=6cm]{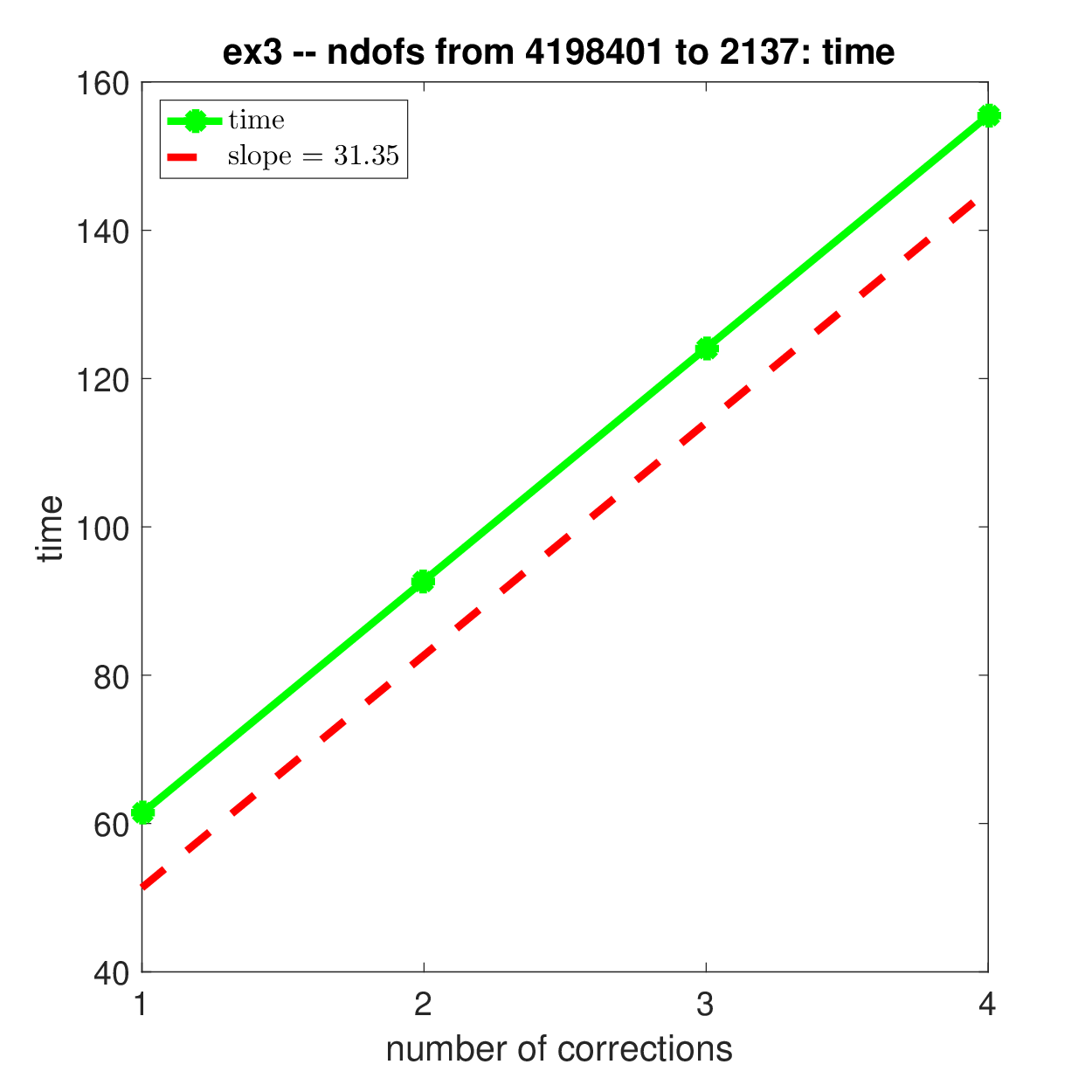}
\caption{Problem (\ref{poisson problem}), example 3 -- The algebraic errors
and CPU time (in second) of the AMG method for the first $13$ eigenvalues
on the uniform refinement mesh.}\label{example 3 figure}
\end{center}
\end{figure}

Table \ref{example 3 table} is presented to show numerical results for various desired eigenvalues.
In order to see the {\it globally uniform convergence} ratio, the dimension on the coarsest
level are increased to $2137$. We also plot the total errors and CPU time for the first $13$
eigenvalues as a function of the number of corrections, respectively, in Figure \ref{example 3 figure}.
We can draw the conclusions that the decreasing of the total error keeps robust
with the number of corrections and the CPU time is still very nice,
growing linearly with the  number of corrections.

\begin{table}[http]
\begin{center}
\begin{tabular}{c|c|c|c|c|c}
\hline
$q$ &       $\lambda_q$     & $\lambda_q / \lambda_{q+1}$ &    ratio   & iter & total error  \\ \hline
1  &    0.019726703793271  &    0.399974375039268        &  0.097224  &  3   &   0.10e-09   \\
2  &    0.049319919035649  &    		                  &            &      &              \\
3  &    0.049320060970441  &    0.624839632041262        &  0.096406  &  3   &   0.66e-09   \\
4  &    0.078932350704643  &    0.800221419880202        &  0.098629  &  4   &   0.11e-09   \\
5  &    0.098638137824978  &    		                  &            &      &              \\
6  &    0.098638194206567  &    0.768992275837269        &  0.097949  &  4   &   0.23e-09   \\
7  &    0.128269421300976  &    		                  &            &      &              \\
8  &    0.128270399246696  &    0.764957776222445        &  0.097399  &  4   &   0.38e-09   \\
9  &    0.167682979680431  &    		                  &            &      &              \\
10  &    0.167683133739996  &    0.944062075515168        &  0.096656  &  4   &   0.57e-09   \\
11  &    0.177618758436507  &    0.900067811260723        &  0.096447  &  4   &   0.68e-09   \\
12  &    0.197339307343651  &    		                  &            &      &              \\
13  &    0.197339360277547  &    0.799924741492187        &  0.096014  &  4   &   0.91e-09   \\
14  &    0.246697407945438  &    		                  &            &      &              \\
15  &    0.246701139291675  &    0.961965508840176        &  0.098786  &  5   &   0.13e-09   \\
16  &    0.256455285584114  &    		                  &            &      &              \\
17  &    0.256455427715290  &    0.896251379863420        &  0.098347  &  5   &   0.16e-09   \\
18  &    0.286142296098190  &    		                  &            &      &              \\
19  &    0.286142698738707  &    0.906120868759646        &  0.097962  &  5   &   0.20e-09   \\
20  &    0.315788664188252  &    0.941206046642901        &  0.097811  &  5   &   0.22e-09   \\
21  &    0.335514912292169  &    		                  &            &      &              \\
22  &    0.335514964638074  &    0.919331113066801        &  0.097356  &  5   &   0.26e-09   \\
23  &    0.364955520235607  &    		                  &            &      &              \\
24  &    0.364955778306734  &    0.924686470880238        &  0.096885  &  5   &   0.31e-09   \\
25  &    0.394680564493737  &    		                  &            &      &              \\
26  &    0.394680661391869  &    0.975464015890347        &  0.104709  &  5   &   0.50e-09   \\
27  &    0.404608119789665  &    		                  &            &      &              \\
28  &    0.404618330771774  &    0.911165080572709        &  0.096076  &  5   &   0.41e-09   \\
29  &    0.444066985663512  &    		                  &            &      &              \\
30  &    0.444068055608404  &    0.900409760718495        &  0.095797  &  5   &   0.47e-09   \\ \hline
\end{tabular}
\caption{Problem (\ref{poisson problem}), example 3 -- Results about the algebraic errors on
unit square with uniform refinement mesh.}\label{example 3 table}
\end{center}
\end{table}

\subsection{Example 4: Poisson eigenvalue problem with discontinuous parameters II}
Here, we still consider eigenvalue problem (\ref{poisson problem}) with coefficient matrix
\begin{equation*}
\left\{
\begin{array}{ccll}
 K &=& \begin{pmatrix}    1  &  0 \\  0 &   1 \end{pmatrix} \ &{\rm in}\ (-1,0) \times ( 0,1) \cup ( 0,1) \times (-1,0),\\
K &=& \begin{pmatrix}   10  &  0 \\  0 &  10 \end{pmatrix} \ &{\rm in}\ [ 0,1] \times [ 0,1] \cup [-1,0] \times [-1,0].
\end{array}
\right.
\end{equation*}
The mesh is generated by uniform refinement beginning with mesh size $h=1$.

\begin{figure}[http]
\begin{center}
\includegraphics[width=7cm,height=6cm]{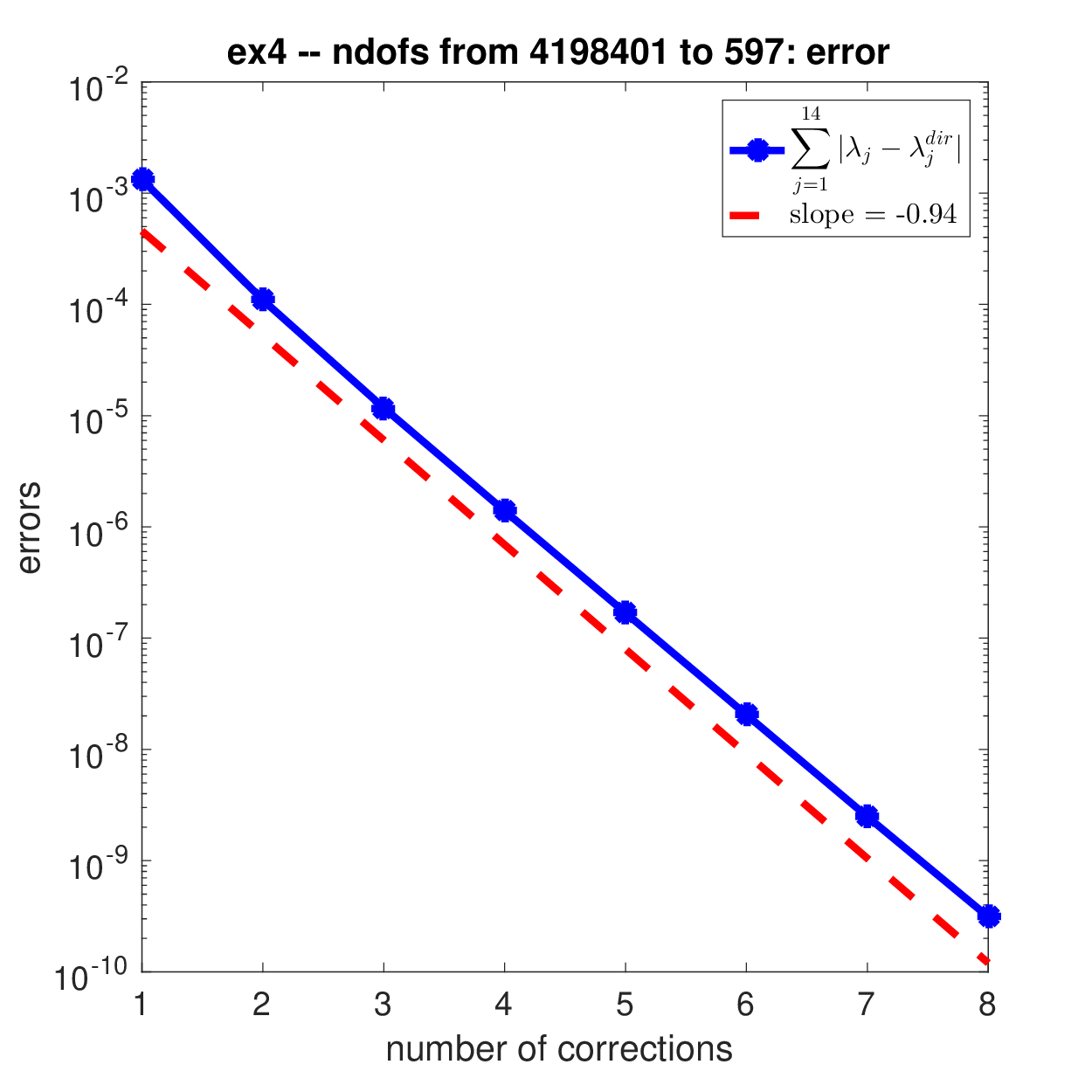}
\includegraphics[width=7cm,height=6cm]{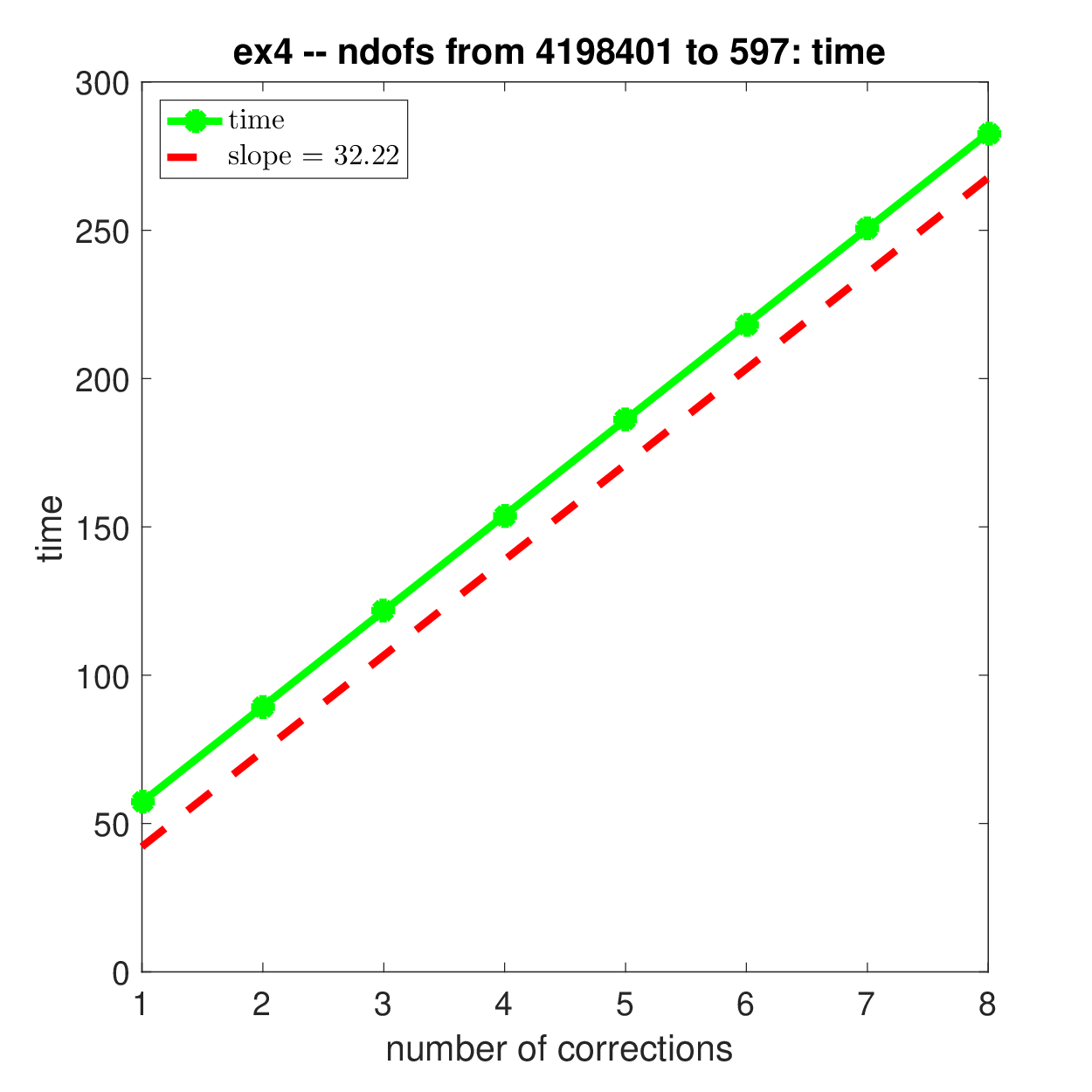}
\caption{Problem (\ref{poisson problem}), example 4 -- The algebraic errors
and CPU time (in second) of the AMG method for the first $14$ eigenvalues
on the uniform refinement mesh.}\label{example 4 figure}
\end{center}
\end{figure}

Different from former three examples, we do not observe similar {\it globally uniform convergence} behavior
even if we increase the coarsest level dimension to $8515$.
In practical computing, we could take at least two simple strategies to improve the performance:
increase the coarsest level dimension, or compute extra eigenpairs.
In Table \ref{example 4 table}, we present numerical experiments of the first 14 eigenvalues taking the two improvement strategies.
Note that the gap between the 14th eigenvalue and the 15th eigenvalue is
$\lambda_{14} / \lambda_{15} = 0.873102340090717$.
Two AMG hierarchies are generated, where dimensions on each level are [4198401, 2095105, 525311, 131589, 33028, 8515, 2165, 597, 176] and
[4198401, 2095105, 525311, 131589, 33028, 8515, 2165, 597], respectively. And three extra eigenpairs,
namely the first 17 eigenvalues and the corresponding eigenvectors are computed as a comparison.
The algebraic errors and CPU time (in seconds) for the first 14 eigenvalues for the third situation,
which takes the shortest time, is given in Figure \ref{example 4 figure}.
Compared with Example 3, we see a big increase in the number of iterations needed to achieve the same tolerance
for different discontinuous parameters, so does the CPU time.

\begin{table}[http]
\begin{center}
\begin{tabular}{c|c|c|c|c|c|c}
\hline
$q$  & actual $q$ & $d_n$  &    ratio   & iter & total error  & time   \\ \hline
14  &   14         & 176  &  0.252941  & 12   &  0.39e-09    & 406.1  \\
14  &   17         & 176  &  0.132108  &  8   &  0.96e-09    & 336.9  \\ \hline
14  &   14         & 597  &  0.113013  &  8   &  0.31e-09    & 282.7  \\
14  &   17         & 597  &  0.113353  &  8   &  0.27e-09    & 359.9  \\ \hline
\end{tabular}
\caption{Problem (\ref{poisson problem}), example 4 -- Results about the algebraic errors
on unit square with uniform refinement mesh.}\label{example 4 table}
\end{center}
\end{table}

\begin{remark}
If the GMG method in \cite{LinXie_2011,Xie_JCP,Xie_IMA} is applied to solve  problem (\ref{poisson problem})
with the above coefficient matrix $K$ defined in this example
and the hierarchy levels are [1050625, 263169, 66049, 16641, 4225, 1089,  289],
we cannot find the convergence behavior, see Table \ref{example 4 table gmg} for details.
The reason should come from that the algebraic multigrid hierarchy has the adaptive coarsening property,
according to the coefficient functions.
\end{remark}

\begin{table}[http]
\begin{center}
\begin{tabular}{c|c|c}
\hline
{correction} & {$\sum_{j=1}^{6}|\lambda_j-\lambda_j^{dir}|$} & ratio \\
\hline
1  &   2.526094e+01  &    -     \\ \hline
2  &   2.565658e+01  &  1.015662\\ \hline
3  &   2.421893e+01  &  0.943965\\ \hline
4  &   2.499157e+01  &  1.031902\\ \hline
5  &   2.417530e+01  &  0.967338\\ \hline
6  &   2.556458e+01  &  1.057467\\ \hline
\end{tabular}
\caption{Problem (\ref{poisson problem}), example 4 -- Results about the number of
the correction and the corresponding algebraic errors
on unit square with uniform refinement mesh for the first 6 eigenvalues.}\label{example 4 table gmg}
\end{center}
\end{table}

In order to show the generality of the AMG method, Delaunay scheme is applied to generate
the unstructured meshes which has no hierarchical structure in the following examples.
Similarly to Example 4, we also test the performance with the two improvement strategies mentioned above.

\subsection{Example 5: Example 1 with Delaunay mesh}

\begin{figure}[http]
\begin{center}
\includegraphics[width=7cm,height=6cm]{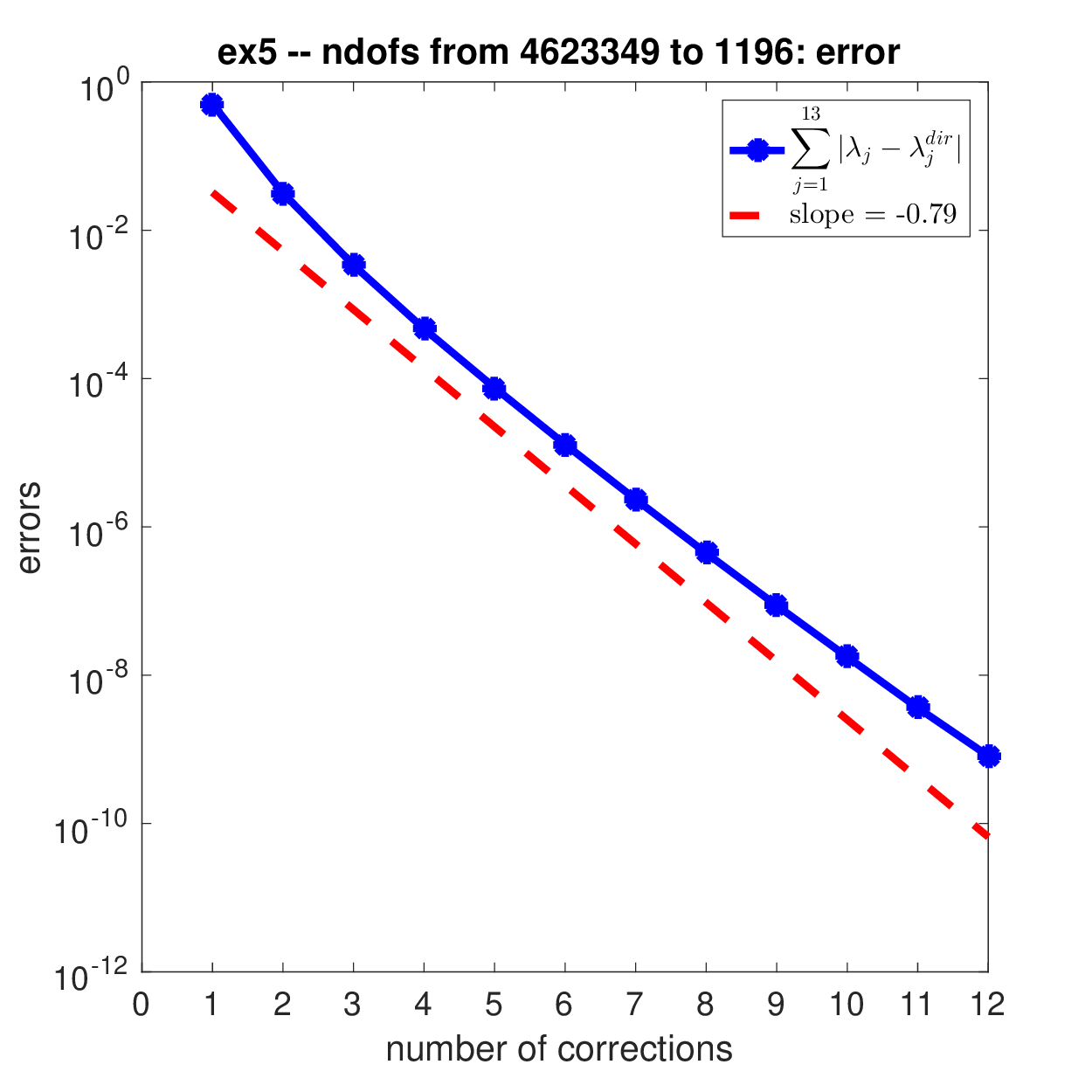}
\includegraphics[width=7cm,height=6cm]{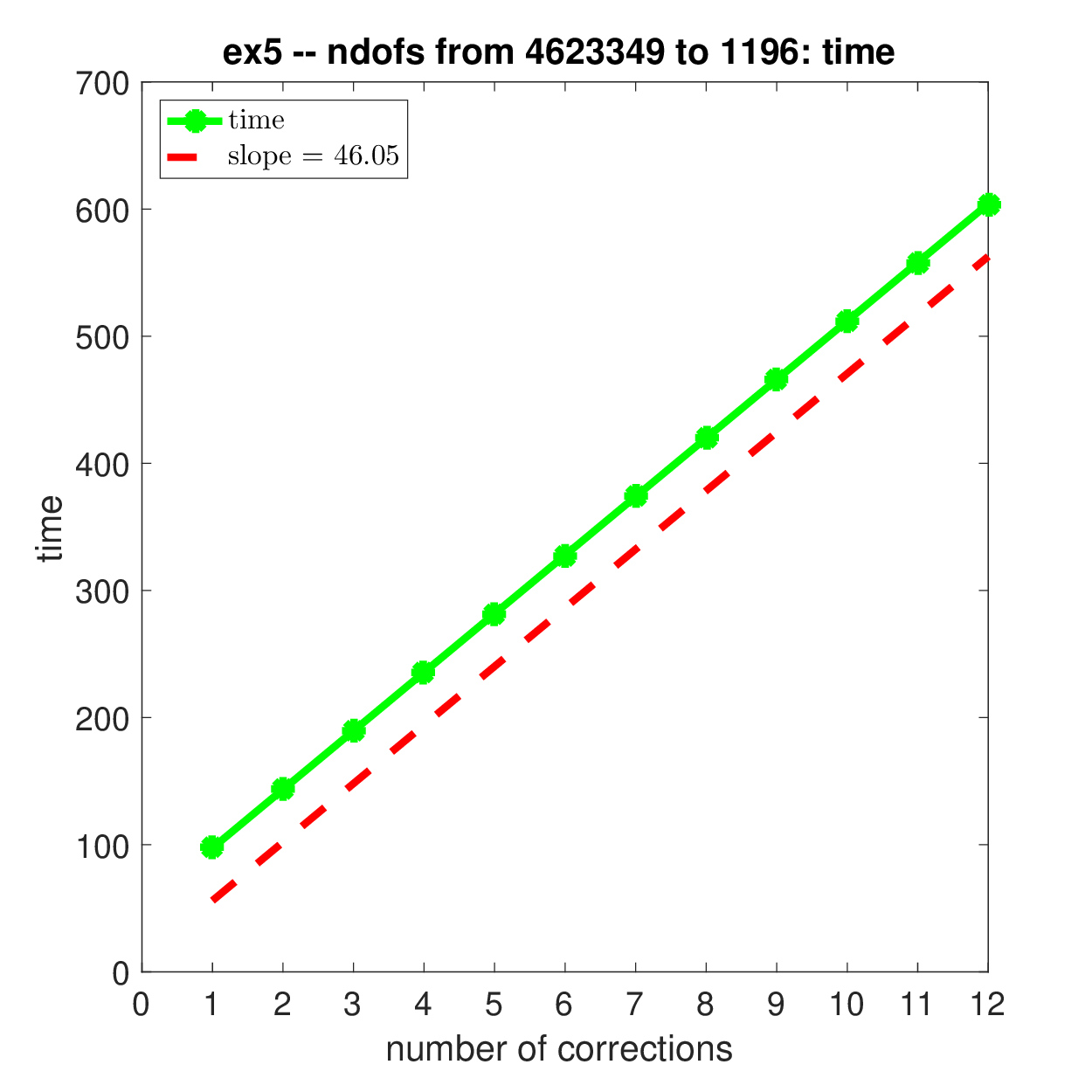}
\caption{Problem (\ref{model problem}), example 5 -- The algebraic errors
and CPU time (in second)  of the AMG method for the first $13$ eigenvalues
on the Delaunay mesh.}\label{example 5 figure}
\end{center}
\end{figure}

Here, we solve the model eigenvalue problem (\ref{model problem}) on the unit square.
Table \ref{example 5 table} presents the numerical results for computing the first 13 eigenvalues.
Note that the gap between the 13th eigenvalue and the 14th eigenvalue is
$\lambda_{13} / \lambda_{14} = 0.799999552462451$.
Dimensions of two AMG hierarchy  on each level are
[4623349, 1630814, 614806, 222577, 76384, 27104, 9743, 3459, 1196] and
[4623349, 1630814, 614806, 222577, 76384, 27104, 9743, 3459], respectively.
For comparison, the first 17 eigenvalues and the corresponding eigenvectors are also computed.
Figure \ref{example 5 figure} reports the algebraic errors and CPU time (in seconds)
for the first 13 eigenvalues with the second situation costing the shortest time in Table \ref{example 5 table}.
\begin{table}[http]
\begin{center}
\begin{tabular}{c|c|c|c|c|c|c}
\hline
$q$ & actual $q$ & $d_n$  &    ratio   & iter & total error  & time   \\ \hline
13  &   13       & 1196   &  0.305067  & 18   &  0.90e-09    & 659.3  \\
13  &   17       & 1196   &  0.158769  & 12   &  0.81e-09    & 604.0  \\ \hline
13  &   13       & 3459   &  0.225724  & 15   &  0.44e-09    & 650.5  \\
13  &   17       & 3459   &  0.151409  & 12   &  0.46e-09    & 723.2  \\ \hline
\end{tabular}
\caption{Problem (\ref{model problem}), example 5 -- Results about the algebraic
errors on unit square with Delaunay mesh.}\label{example 5 table}
\end{center}
\end{table}

\subsection{Example 6: Example 2 with Delaunay mesh}
\begin{figure}[http]
\begin{center}
\includegraphics[width=7cm,height=6cm]{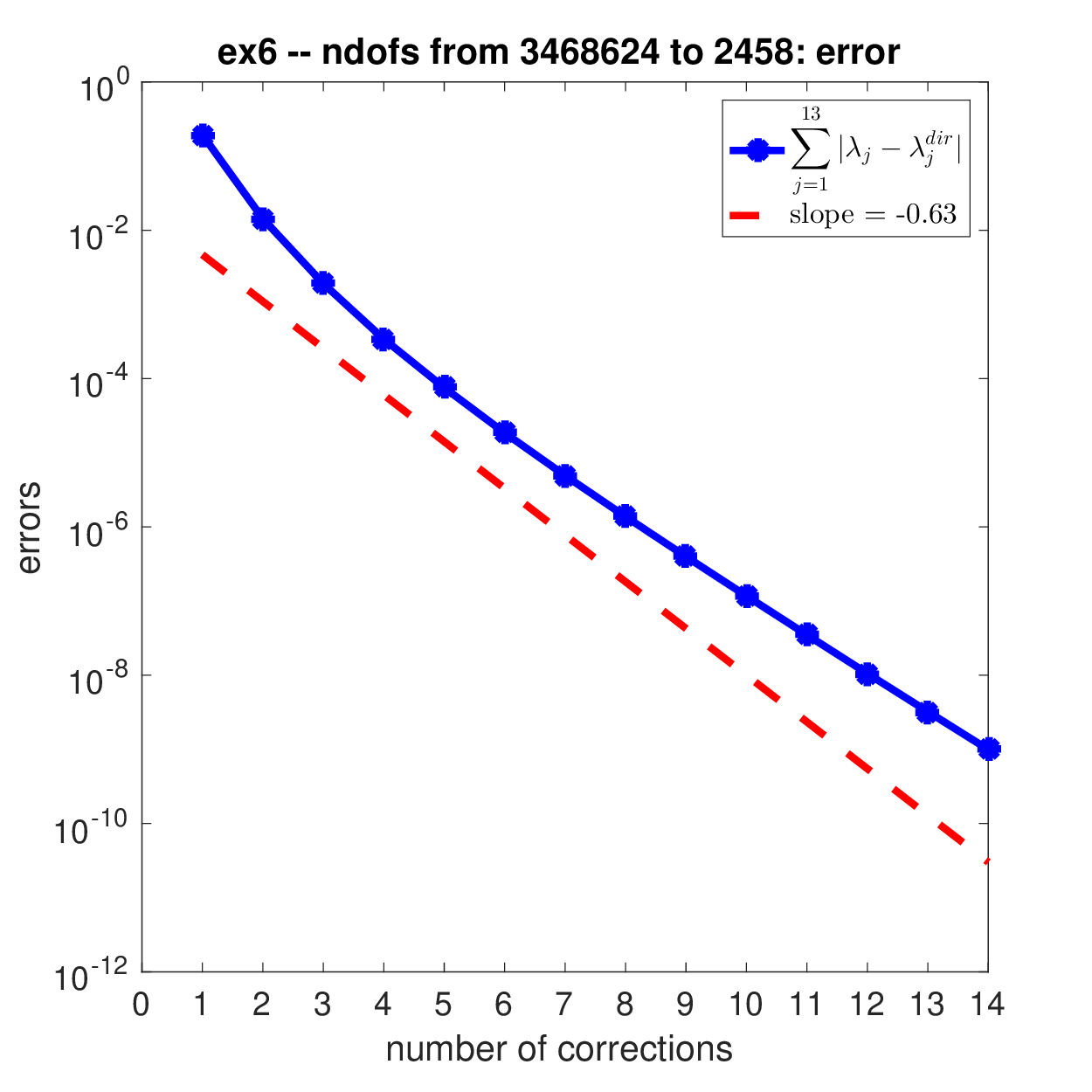}
\includegraphics[width=7cm,height=6cm]{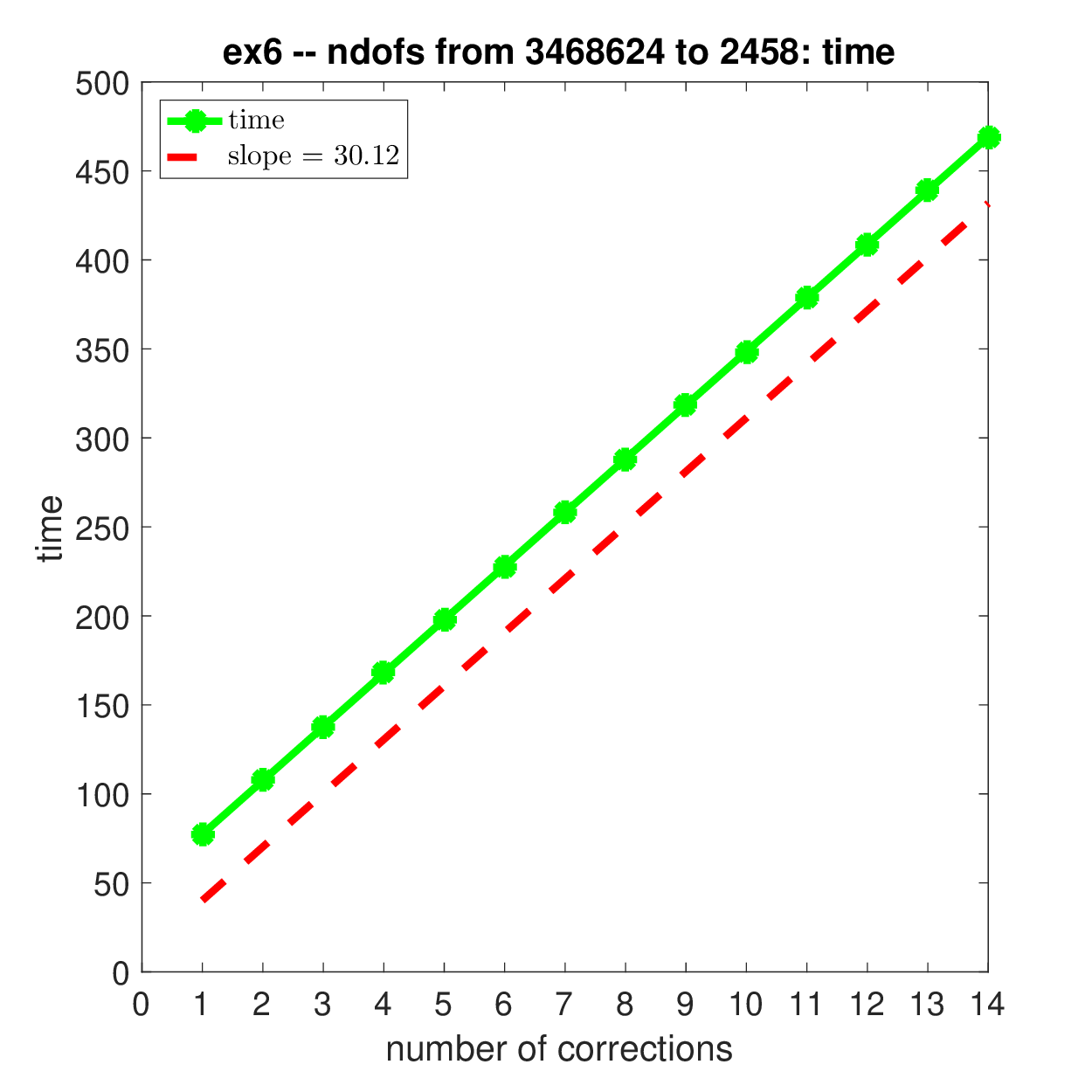}
\caption{Problem (\ref{model problem}), example 6 -- The algebraic errors  and CPU time (in second)
of the AMG method for the first $13$ eigenvalues
on the Delaunay mesh.}\label{example 6 figure}
\end{center}
\end{figure}

We solve the model eigenvalue problem (\ref{model problem})
on the L-shape domain $\Omega=[-1,1]\times [-1,1]\backslash(0,1)\times(-1,0)$ with the mesh generated by Delaunay method.
In Table \ref{example 6 table}, we present the numerical results for computing the first 13 eigenvalues.
Note that the gap between the 13th eigenvalue and the 14th eigenvalue is
$\lambda_{13} / \lambda_{14} = 0.906477598649387$. Dimensions on each level are
[3468624, 1219797, 453107, 165008, 55740, 19564, 6958, 2458, 865] and
[3468624, 1219797, 453107, 165008, 55740, 19564, 6958, 2458], respectively.
Similarly to the former example, various extra eigenvalues, from the first 14 to 17 eigenvalues,
are computed for comparisons. In this example, the second situation takes the shortest time in Table \ref{example 6 table}.
Figure \ref{example 6 figure} shows the correspond convergence behavior and CPU time (in seconds).

\begin{table}[http]
\begin{center}
\begin{tabular}{c|c|c|c|c|c|c}
\hline
$q$  & actual $q$ & $d_n$  &    ratio   & iter & total error  & time   \\ \hline
13  &   13       & 865    &  0.540191  & 20   &  0.20e-05    & 525.8  \\
13  &   14       & 865    &  0.308052  & 18   &  0.42e-09    & 517.2  \\
13  &   15       & 865    &  0.247463  & 15   &  0.66e-09    & 476.7  \\
13  &   16       & 865    &  0.219253  & 14   &  0.53e-09    & 481.3  \\
13  &   17       & 865    &  0.214510  & 14   &  0.40e-09    & 510.0  \\ \hline
13  &   13       & 2458   &  0.468039  & 20   &  0.11e-06    & 597.9  \\
13  &   14       & 2458   &  0.230357  & 14   &  0.99e-09    & 469.0  \\
13  &   15       & 2458   &  0.182574  & 13   &  0.26e-09    & 469.2  \\
13  &   16       & 2458   &  0.159574  & 12   &  0.31e-09    & 471.8  \\
13  &   17       & 2458   &  0.155989  & 12   &  0.24e-09    & 507.2  \\ \hline
\end{tabular}
\caption{Problem (\ref{model problem}), example 6 -- Results about the algebraic errors on
an L-shape domain with Delaunay mesh.}\label{example 6 table}
\end{center}
\end{table}

\section{Conclusions}
A type of AMG method to solve algebraic eigenvalue problems
arising from the discretization of partial differential equations is developed. Paired with
the multilevel correction method and the AMG method for linear equations, the resulting algorithm
for eigenvalue problems needs almost the optimally computational work and the least memory.
The efficiency of the proposed AMG method is exhibited by six numerical examples. The first three examples are presented to show that the AMG method has {\it globally uniform convergence} rate under some good situations.
Another three challenging examples are included to discuss strategies to improve the convergence under not-so-good situations.
In our implementation, the choices of interpolation types (direct or standard), pre-
and postsmoothing operators (CG, GS, etc.), linear solvers (AMG V-cycle, CG, PCG, etc.), and parameters $p_k$ and $d_n$, etc.,
are tested.

In order to develop a robust and efficient AMG eigensolver for eigenvalue problems, different type of coarsening and interpolation
strategies should be tested for different type of matrices, other types of eigenvalue solution schemes, such as
shift strategy, Lanzcos method, restarting techniques, polynomial acceleration, should be embedded in the AMG method here.
Furthermore, for large scale eigenvalue problems, parallelization implementation should also be considered and implemented.
These will be our further work.

\section*{Acknowledgments}
We would like to express our sincerely thanks to Prof. Chensong Zhang for providing the algebraic multigrid software.

The research has been supported by Science Challenge Project (No. TZ2016002),
National Natural Science Foundations of China (NSFC 11771434, 91330202, 11371026, 11001259),
the National Center for Mathematics and Interdisciplinary Science, CAS.

\end{document}